\documentclass[12pt]{amsart}
\usepackage{latexsym}
\usepackage{amssymb}
\usepackage{amscd}
\renewcommand\a{\alpha}
\renewcommand\b{\beta}

\newcommand\la{\lambda}
\newcommand\z{\zeta}

\renewcommand\th{\theta}

\newcommand\f{\phi}

\newcommand\p{\psi}
\renewcommand\t{\tau}

\newcommand\F{\varPhi}

\newcommand\vG{\varGamma}

\newcommand\Ql{\bar{\mathbf Q}_l}
\newcommand\BQ{\mathbf Q}
\newcommand\BF{\mathbf F}
\newcommand\BC{\mathbf C}
\newcommand\BR{\mathbf R}
\newcommand\BZ{\mathbf Z}

\newcommand\CB{\mathcal{B}}

\newcommand\CH{\mathcal{H}}
\newcommand\CI{\mathcal{I}}
\newcommand\CJ{\mathcal{J}}

\newcommand\CL{\mathcal{L}}

\newcommand\CP{\mathcal{P}}

\newcommand\CZ{ \mathcal{Z}}

\newcommand\FS{\mathfrak S}

\newcommand\Fp{\mathfrak p}
\newcommand\iv{^{-1}}

\newcommand\wt{\widetilde}

\newcommand\IC{\operatorname{IC}}

\newcommand\Ind{\operatorname{Ind}}

\newcommand\Tr{\operatorname{Tr}\,}

\newcommand\ad{\operatorname{ad}}

\newcommand\der{\operatorname{_{\!der}}}
\newcommand\reg{_{\operatorname{reg}}}

\newcommand\lp{\operatorname{\!\langle\!}}
\newcommand\rp{\operatorname{\!\rangle\!}}

\newcommand\odd{\operatorname{odd}}
\newcommand\Gal{\operatorname{Gal}}
\newcommand\ch{\operatorname{ch}}

\newcommand\da{\dot a}
\newcommand{\dda}{\ddot a}

\newcommand{\isom}{\,\raise2pt\hbox{$\underrightarrow{\sim}$}\,}

\newcommand{\gind}{\vG\!\operatorname{-Ind}}
\numberwithin{equation}{section}

\newtheorem{thm}{Theorem}[section]
\newtheorem{lem}[thm]{Lemma}
\newtheorem{cor}[thm]{Corollary}
\newtheorem{prop}[thm]{Proposition}

\def \para#1{\par\medskip\textbf{#1}
              \addtocounter{thm}{1}}

\def \remark#1{\par\medskip\noindent
                \textbf{Remark #1}
                \addtocounter{thm}{1}}

\def \remarks#1{\par\medskip\noindent
                \textbf{Remarks #1}
                \addtocounter{thm}{1}}

\begin{document}
\setlength{\baselineskip}{4.9mm}
\setlength{\abovedisplayskip}{4.5mm}
\setlength{\belowdisplayskip}{4.5mm}
\renewcommand{\theenumi}{\roman{enumi}}
\renewcommand{\labelenumi}{(\theenumi)}
\renewcommand{\thefootnote}{\fnsymbol{footnote}}
\renewcommand{\thefootnote}{\fnsymbol{footnote}}
\parindent=20pt
\medskip
\begin{center}
{\bf A variant of the induction theorem \\ for Springer representations} \\
\vspace{1cm}
Toshiaki Shoji\footnote{
A part of this work was done during the author's stay in 
CIB in Lausanne in June 2005.  The author is grateful for 
their hospitality.} 
\\ 
\vspace{0.5cm}
Graduate School of Mathematics \\
Nagoya University  \\
Chikusa-ku, Nagoya 464-8602,  Japan
\end{center}
\title{}
\maketitle
\begin{abstract}
Let $G$ be a simple algebraic group over $\BC$ with the Weyl group
$W$.  For a unipotent element $u \in G$, let $\CB_u$ be the
variety of Borel subgroups of $G$ containing $u$.  
Let $L$ be a Levi subgroup of a parabolic subgroup of $G$ with the 
Weyl subgroup $W_L$ of $W$.
Assume that $u \in L$ and let $\CB_u^L$ be a similar variety as $\CB_u$ 
for $L$. 
For a certain choice of $L$, $u \in L$ 
and $e \ge 1$, we describe the $W$-modules 
$\bigoplus_{n \equiv k \mod e}H^{2n}(\CB_u)$  for 
$k = 0, \dots, e-1$, in terms of the $W_L$-module $H^*(\CB_u^L)$
with some additional data, which
is a refinement of the induction theorem due to Lusztig.   
As an application, we give an explicit formula for the values of Green
functions at root of unity, in the case where $u$ is a regular
unipotent element in $L$.
\end{abstract}
\pagestyle{myheadings}
\markboth{SHOJI}{SPRINGER REPRESENTATIONS}

\bigskip 
\medskip
\addtocounter{section}{-1}
\section{Introduction}
Let $G$ be a connected reductive group over an algebraically closed
field $k$, and $W$ the Weyl group of
$G$.  For a unipotent element $u \in G$, let $\CB_u$ be the variety of
Borel subgroups containing $u$.   According to Springer [Sp2], 
Lusztig [L1], $W$ acts naturally on the $l$-adic cohomology group 
$H^n(\CB_u) = H^n(\CB_u,\Ql)$, the so-called Springer representations
of $W$.  
Assume that $k = \BC$, or the characteristic $p$
of $k$ is good.  Then it is known that $H^{\odd}(\CB_u) = 0$.  
We consider the graded $W$-module 
$H^*(\CB_u) = \bigoplus_{n \ge 0}H^{2n}(\CB_u)$.
Let $L$ be a Levi subgroup of
a parabolic subgroup of $G$. Let $W_L$ be the Weyl group of $L$, which
is naturally a subgroup of $W$. If $u \in L$, the variety 
$\CB_u^L$ is defined by replacing $G$ by $L$, and we have a graded 
$W_L$-module $H^*(\CB_u^L)$.  
\par  
Lusztig proved in [L3] an induction theorem for Springer
representations, which describes the $W$-module structure 
of $H^*(\CB_u)$ in terms of the $W_L$-module
structure of $H^*(\CB_u^L)$, in the case where $u \in L$.
However in this theorem, the information on 
the graded $W$-module structure is eliminated.  
In this paper, we try to recover partly the graded $W$-module structure,
i.e., for a fixed positive integer $e$, 
we consider the $W$-modules 
$V_{e,k} = \bigoplus_{n \equiv k \mod e}H^{2n}(\CB_u)$     
for $k = 0, \dots, e-1$.  Let $G$ be a simple group
modulo center defined over $\BC$.  We show,  under a certain choice of 
$L$, $u$ and $e$, that the $W$-module $V_{e,k}$ can be described in terms 
of the graded $W_L$-module $H^*(\CB_u^L)$ with some additional data.  
In particular, 
we see that $\dim V_{e,k}$ is independent of the choice of $k$. 
\par
In the case where $u = 1$, $H^*(\CB_u)$ is isomorphic, as a graded
$W$-module,  to the coinvariant algebra of $W$.
In this case $V_{e,k}$ has been studied by many authors, by Stembridge [St]
for $e$ corresponding to the regular elements in $W$, by Morita and 
Nakajima [MN1] for $W = \FS_n$ with $e$ such that $1 \le e \le n$, 
and by Bonnaf\'e, Lehrer and Michel [BLM] for complex reflection groups
$W$ in the most general framework.  Our result partly covers the
result of [BLM].  For general $u \ne 1$, Morita and Nakajima [MN2] 
considered certain types of unipotent elements for $G = GL_n$, which is 
a special case of ours. 
\par
The proof of the induction theorem in [L3] is done by passing to the finite
field $\BF_q$, and using a certain specialization argument $q \mapsto 1$ 
together with the properties of Deligne-Lusztig's virtual character $R_T(1)$.     
Our argument is a variant of that in [L3].  We use a
specialization $q \mapsto \z$, where $\z$ is a primitive $e$-th root
of unity.
Thus our argument is closely related to the values of Green functions
at root of unity.  In the case where $u$ is a regular unipotent
element in $L$, we obtain an explicit formula for such values, which 
is regarded as a generalization of the result by 
Lascoux, Leclerc and Thibon [LLT] for the case of 
Green polynomials of $GL_n$.
\par\medskip
\section{The statement of the main result}
\para{1.1.}
Let $k$ be an algebraic closure of a finite field with ch$(k) = p > 0$
or the complex number field $\BC$.  Let $G$ be a connected reductive 
group $G$ over $k$.  
Let $\CB$ be the variety of Borel subgroups of $G$, and $W$ 
the Weyl group of $G$. 
For any $g \in G$, put 
$\CB_g = \{ B' \in \CB \mid g \in B'\}$.   
We consider the Springer representations of $W$ on 
$H^n(\CB_g, \Ql)$ (or on $H^n(\CB_g, \BC)$ in the case where $k = \BC$).
\par
Let $L$ be a Levi subgroup of a parabolic subgroup $P$ of $G$. 
The  Weyl group $W_L$ of $L$ is naturally identified with a subgroup
of $W$. 
Let $\CB^L$ be the variety of Borel subgroups of $L$. 
For a unipotent element $u \in L$, we consider
$\CB_u^L = \{ B' \in \CB^L \mid u  \in B'\}$. Thus we have a $W_L$-module 
$H^n(\CB^L_u,\Ql)$, and a $W$-module $H^n(\CB_u, \Ql)$.
The induction theorem for Springer representations asserts that
\begin{equation*}
\tag{1.1.1}
\sum_{n \ge 0}(-1)^nH^n(\CB_u, \Ql) = 
   \Ind_{W_L}^W\bigl(\sum_{n \ge 0}(-1)^nH^n(\CB^L_u, \Ql)\bigr)
\end{equation*}
as virtual $W$-modules.
\remark{1.2.}
The induction theorem was stated in [AL], with a brief indication 
of the proof, in the case where $k = \BC$, and was proved in [L3]
for any $k$.  Note that if $p$ is good, the unipotent classes in 
$G$ are parametrized in the same way as the case of $k = \BC$,
independent of $p$.  Moreover in that case, it is known that 
$H^n(\CB_u,\Ql) = 0$ for odd $n$.  Then the algorithm of computing 
Green functions implies that the $W$-module structure of 
$H^n(\CB_u,\Ql)$ is independent of $p$.  Thus by a general principle 
$H^n(\CB_u,\Ql)$ is isomorphic to the $W$-module 
$H^n(\CB_{u'},\BC)$, where $u', \CB_{u'}$ are the corresponding 
objects in the algebraic group $G_{\BC}$ over $\BC$. 
In what follows, we express $H^n(\CB_u, \Ql)$ or 
$H^n(\CB_{u'}, \BC)$ by $H^n(\CB_u)$ by abbreviation.
\para{1.3.}
Assume that $k = \BC$.
We consider the following variant of the induction theorem. 
Let $\vG$ be a cyclic group of order $e$ generated by $a$.  
Let $\z$ be a primitive $e$-th root of unity in $\BC$.
Let $V = \bigoplus_{n\ge 0}V_n$ be a graded $W$-module.  Then 
$V$ turns out to be a $\vG \times W$-module by defining the 
action of $\vG$ on $V$ by $ax = \z^nx$ for $x \in V_n$.
We denote by $V^{(\z)}$ the thus obtained $\vG \times W$-module $V$. 
\par
For $u \in L$, we consider 
the graded $W_L$-module 
$H^*(\CB^L_u) = \bigoplus_{n \ge 0}H^{2n}(\CB^L_u)$, where 
the degree $n$ part is given by $H^{2n}(\CB^L_u)$, 
and similarly we consider the graded $W$-module 
$H^*(\CB_u) = \bigoplus_{n \ge 0}H^{2n}(\CB_u)$. 
Let $\vG$ be as before.  
We choose $\vG$ such that $\vG \subset N_W(W_L)$, and 
consider the semidirect product $\wt W_L = \vG \ltimes W_L$. 
We assume that the $W_L$-module $H^n(\CB_u^L)$ can be extended to
a $\wt W_L$-module for each $n$.  (In the case where $a \in Z_W(W_L)$, 
we have $\wt W_L = \vG \times W_L$.  In this case, one can choose
a trivial extension to $\wt W_L$, i.e., we may asssume that 
$\vG (\subset \wt W_L)$ acts trivially on $H^*(\CB_u^L)$.)  
Then one can define a
$\vG \times \wt W_L$-module $H^*(\CB_u^L)$ as above, replacing $W_L$ by 
$\wt W_L$, which we denote
by $H^*(\CB_u^L)^{(\z)}$.
(When we need to distinguish the group $\vG$ as the first factor of 
$\vG \times \wt W_L$ from the subgroup of $\wt W_L$, 
we write the latter as $\vG_0$.)
$\vG \times W$-module $H^*(\CB_u)^{(\z)}$ is defined as before.
Put $V^{(\z)} = H^*(\CB^L_u)^{(\z)}$, and let $V^{(\z)}_n$ be the degree 
$n$-part of $V^{(\z)}$.  Let us consider the induced 
$W$-module 
\begin{equation*}
\Ind_{W_L}^WV^{(\z)} = \bigoplus_{w \in W/W_L}w\otimes V^{(\z)}.
\end{equation*}
Then $\Ind_{W_L}^WV^{(\z)}$ turns out to be a 
$\vG \times W$-module by defining the action of $\vG$ by 
$b(w\otimes x) = \z^n(wb\iv \otimes bx)$ for $b \in \vG_0, x \in V_n^{(\z)}$, 
which we denote by $\gind_{W_L}^WV^{(\z)}$.
\para{1.4.}
In the remainder of this paper, we assume that $G$ is simple modulo
center. 
Let $T \subset B$ be a pair of maximal torus and a Borel subgroup of 
$G$. Put $W = N_G(T)/T$.  
Let $L$ be a Levi subgroup of a parabolic subgroup $P$ of $G$
containing $B$ such that $L \supset T$. We have $W_L = N_L(T)/T$.
Let $\F \subset X(T)$ be a root system for $G$ with respect to $T$, 
with a simple root system $\Pi$ (with respect to $B$), where 
$X(T)$ is the character group of $T$.
We denote by $\F_L$ the sub system of $\F$ corresponding to 
$L$ with the simple root system  $\Pi_L \subset \Pi$.
Let $\Pi'$ be the set of simple roots which are orthogonal to $\Pi_L$
with respect to the standard inner product on 
$V = \BR\otimes_{\BZ}X(T)$.
We denote by $L'$ the Levi subgroup containing $T$ 
corresponding to $\Pi'$.  Let $W_{L'} = N_{L'}(T)/T$ be the Weyl group
of $L'$.  Then we have $W \supset W_L \times W_{L'}$, and so
$W_{L'} \subset N_W(W_L)$.
\par
We recall here the notion of regular elements of reflection groups 
due to Springer [Sp1]. Let $W$ be a reflection group in $GL(V)$.
A vector $v \in V$ is called regular if $v$ is not contained in any 
reflecting hyperplane in $V$.  An element $a \in W$ is called regular 
if $a$ has an eigenvector $v$ which is a regular element in $V$. 
If $av = \z v$, with $\z$ a primitive $e$-th root of unity, 
then the order of $a$ is equal to $e$ ([Sp1, 4.2]).  In particular,
if $a$ is regular of order $e$, there exists an eigenvalue $\z$ 
which is a primitive $e$-th root of unity.
\par
The regular elements $a \in W$ 
in the case of classical groups are given as
follows (cf. [Sp1]). 
\par
{\bf Type $A_{n-1}$}.  In this case $W= \FS_n$ and there are two types
of regular elements.
\par
(a) \ $e$ is a divisor of $n$, and $a$ is an $n/e$-product
of (disjoint) $e$-cycles in $\FS_n$.
\par
(b) \ $e$ is a divisor of $n-1$, and $a$ is an $(n-1)/e$-product 
of $e$-cycles in $\FS_n$
\par
{\bf Type $B_n$}.  There are two types of regular elements.
\par
(a) \ $e$ is an odd divisor of $n$, and $a$ is an $n/e$-product of
positive cycles of length $e$.  
\par
(b) \ $e$ is an even divisor of
$2n$, and $a$ is a $2n/e$-product of negative cycles of length $e/2$.
\par
{\bf Type $D_n$}. 
In this case there are 4 types of regular elements.
\par
(a) \ $e$ is an odd divisor of $n$, and $a$ is a product of positive
cycles of length $e$.
\par
(b) \ $e$ is an odd divisor of $n-1$, and $a$ is a product of positive
cycle of length 1 and $(n-1)/e$ positive cycles of length $e$.
\par
(c) \ $n$ is even, and $e$ is an even divisor of $n$. $a$ is a product
of negative cycles of length $e/2$.
\par
(d) \ $e$ is an even divisor of $2n-2$, and $a$ is a product of 
$(n-1)/e$ negative cycles of length $e/2$ and one cycle of length 1,
which is positive or negative according as $(2n-2)/e$ is even or odd.
\par  
Regular elements in the exceptional Weyl groups are listed in [Sp1].
\par
Returning to the original setting, we consider the subgroups 
$W_L, W_{L'}$ of $W$.  Let $V'$ be the subspace of $V$ generated by 
$\Pi_{L'}$. $W_{L'}$ is realized as a reflection group on $V'$.  
Assume that $a$ is a regular element of $W_{L'}$ of order $e$. Let 
$\z$ be a primitive $e$-th root of unity, and $V(a, \z)$ the
eigensubspace of $a$ in $V$ with eigenvalue $\z$. 
Since $a$ is regular, $V(a,\z)$ is not contained in any reflecting 
hyperplane $H_{\a}$ for $\a \in \F_{L'}$.  We say that $a$ is 
$L$-regular if $V(a,\z)$ is not contained in any $H_{\a}$
for $\a \in \F - \F_{L}$.
If $L$ is the torus $T$, all the regular elements are $L$-regular. 
But if $L \ne T$, regular elements are not necessarily $L$-regular. 
For example, if $L$ is not simple modulo center, 
regular elements in $W_{L'}$ are not $L$-regular in many cases. 
In the case where $L$ is simple modulo center, 
$L$-regular elements are classified as follows.
\begin{lem}  
Assume that $L$ is simple modulo center.  
\begin{enumerate}
\item
If $W$ is of type $A_n, B_n, D_n$, take $L$ such that  
$W_L$ is of the same type as $W$ of rank $m$, and $W_{L'}$ is of type
$A_{n-m-1}$. Then a regular element of $W_{L'}$ of type
\rm{(a)} in 1.4 is $L$-regular.   
\item
If $W$ is of type $G_2, F_4$ or $E_8$, there does not exist
$L$-regular elements for any $L \ne T$.
\item
Assume that $W$ is of type $E_6$ or $E_7$.  Let 
$\Pi = \{ \a_1, \dots, \a_7\}$ (resp. $\{\a_1, \dots, \a_6\}$) 
be the set of simple roots in $E_7$ (resp. in $E_6$) as in the figure. 
Take $\Pi_L = \{ \a_k, \a_{k+1}, \dots, \a_7\}$ 
(resp. $\{\a_k, \a_{k+1}, \dots, \a_6\}$)
for $k \ge 3$.
Then $W_{L'}$ is of type $A_j$ or of type $A_j + A_1$ for some $j$ 
except the case where $W$ is of type $E_7$ and $\Pi_L = \{ \a_7\}$, 
in which case $\Pi_{L'}$ is of type $D_5$.
In the former case, we choose $a$ a regular element of type 
\rm{(a)} for type $A$, and in the latter case, we choose 
$a$ a regular element of type \rm{(a)} for type $D$ in 1.4, respectively.
Then $a$ is $L$-regular.
\begin{figure}[h]
\setlength{\unitlength}{1mm}
\begin{picture}(80,24)(-20,0)
\multiput(-10,1)(13,0){6}{\thicklines\circle{2}}
\multiput(-9.0,1)(13,0){5}{\thicklines\line(1,0){11}}
\put(16, 12){\thicklines\circle{2}}
\put(16, 2){\thicklines\line(0,1){9}}

\put(-40,1){$E_{7}\quad (E_6)$}

\put(14, 17){$\a_3$}
\put(-12, -4){$\a_1$}
\put(1, -4){$\a_2$}
\put(14,-4){$\a_4$}
\put(27,-4){$\a_5$}
\put(40, -4){$\a_6$}
\put(53, -4){$\a_7$}

\end{picture}
\end{figure}
\end{enumerate}
\end{lem}
\begin{proof}
If there exists $\b \in \F - \F_L$  such that $\b$ is orthogonal
to $V_{L'}$, then any regular element in $W_{L'}$ cannot be 
$L$-regular.  By direct inspections, one can find such $\b$ 
unless $L$ is the type given in (i), (iii) of the lemma. 
Assume that $L$ is
as in the lemma, and let $a$ be a regular element in $W_{L'}$.  
If $W_{L'}$ is of type $A_j$ or type $A_j + A_1$, 
then a regular vector $v \in V'$ can be written explicitly, and
one can check the $L$-regularity by direct inspections. If 
$W_{L'}$ is of type $D_5$ (in the case where $W$ is of type $E_7$),
$a$ must be of type (a) (otherwise it is easy to see that 
$a$ is not $L$-regular).  But this element is nothing but the regular
element in $A_4$, and the checking is reduced to the previous case.    
The details are omitted.
\end{proof}

\para{1.6.}
In what follows we consider a specific cyclic group 
$\vG \in N_W(W_L)$, and $u \in L$ according to the following two cases. 
\par\medskip\noindent
Case (a): $W_{L'} \ne \{ 1\}$.
\par
In this case, we assume that $L$ is simple modulo center. We
choose an $L$-regular element $a \in W_{L'}$, and put 
$\vG = \lp a\rp$.  Let $e$ be the order of $\vG$.
Thus $\vG \subset W_{L'}$ and we have $\vG \times W_L \subset W$.
We take any unipotent element $u \in L$.
\par\medskip\noindent
Case (b): $W_{L'} = \{1\}$.
\par
In this case,  we assume that $L$ is of type 
$X_0 + e(A_{n_1-1} + \cdots + A_{n_r-1})$ 
with $X_0$ irreducible.  
We further assume that any $\b \in \F - \F_L$ is not orthogonal
to the root system $e(A_{n_1-1} + \cdots + A_{n_r-1})$.
(Note: since $W_{L'} = \{1\}$, any irreducible
component of the Dynkin diagram corresponding to $\Pi - \Pi_L$
consists of 1 or 2 nodes.  The latter condition is satisfied 
for type $B_n$ if all the irreducible components consist of one node,
and for type $A_n, D_n$ if the number of irreducible components 
having two nodes is at most 1.)    
\par
We choose $a \in W$ so that $a$ permutes each component 
$A_{n_i-1}$ in a cyclic way, and acts trivially on $X_0$.
Thus $a \in \FS_{en_1}\times \cdots\times \FS_{en_r}$, and $a$
is a product of disjoint cycles of length $e$.  In particular, 
$\vG = \lp a\rp \subset N_W(W_L)$, and the subgroup of $W$ generated
by $\vG$ and $W_L$ coincides with the semidirect product 
$\vG\ltimes W_L$. 
Now $L$ is isogenic to $G_0\times G_1 \times \cdots\times G_r$ modulo
center, where $G_0$ is of type $X_0$, and 
$G_i \simeq GL_{n_i} \times \cdots \times GL_{n_i}$ ($e$-factors).
We choose a unipotent element $u \in L$ so that 
$u$ corresponds to $(u_0, u_1, \dots, u_r)$, where $u_0 \in G_0$ is arbitrary, 
and  $u_i$ is a diagonal element in $G_i$, i.e., 
$u_i = (v_i, \dots, v_i)$ with $v_i \in GL_{n_i}$ for $i = 1, \dots, r$.
\par\medskip
We can state our main theorem, whose proof will be given 
in the next section.  
\begin{thm} 
Assume that $G$ is defined over $\BC$. Let $L$ be a Levi subgroup in $G$. 
Assume that a cyclic subgroup $\vG $ of order $e$ 
in $N_W(W_L)$ and $u \in L$ are given as in 1.4.  
Put $\wt W_L = \vG \ltimes W_L$.  Then the followings hold. 
\begin{enumerate}
\item
$W_L$-module $H^*(\CB_u^L)$ can be extended to a $\wt W_L$-module
so that $\vG \times \wt W_L$-module $H^*(\CB_u^L)^{(\z')}$ is defined
for any $e$-th root of unity $\z'$.
\item
There exists a primitive $e$-th root of
unity $\z$ such that 
\begin{equation*}
\tag{1.7.1}
\gind_{W_L}^W\bigl(H^*(\CB^L_u)^{(\z)}\bigr) \simeq H^*(\CB_u)^{(\z)}
\end{equation*} 
as $\vG \times W$-modules.
\end{enumerate}
\end{thm}
\remarks{1.8.} (i)  The extension of $W_L$-module $H^*(\CB_u^L)$
to $\wt W_L$-module is not unique.  The theorem aaserts that 
the statement (ii) holds for some choice of extension.  
\par
(ii) \ The theorem asserts that (1.7.1) holds for
some choice of primitive $e$-th root of unity $\z$, but then it 
holds for any choice of primitive root of unity $\z'$.
In fact, we can write $\z' = \z^j$ for some $j$ prime to $e$, and 
we have an automorphism $\t$ on $\vG$ such that $\t(a) = a^j$.
It follows from (1.7.1) that we have an isomorphism of $\vG \times W$
modules, where the action of $\vG$ is twisted by $\t$.  
It is easy to check that the twisted $\vG \times W$-module 
$\gind_{W_L}^W\bigl(H^*(\CB^L_u)^{(\z)}\bigr)$ is isomorphic to
$\gind_{W_L}^W\bigl(H^*(\CB^L_u)^{(\z')}\bigr)$, and similarly 
the twisted $H^*(\CB_u)^{(\z)}$ is isomorphic to $H^*(\CB_u)^{(\z')}$.   
Thus (1.7.1) holds also for $\z'$.
\par\bigskip
\section{Proof of Theorem 1.7}
\para{2.1.}  
In the case where $e = 1$, Theorem 1.7 is nothing but the original
induction theorem. So we assume that $e \ge 2$ in what follows. 
Since the structure of the $W$-module $H^n(\CB_u)$ is independent
of $p$ provided that $p$ is a good prime,  
it is enough to show the corresponding formula for 
an appropriate $p$.   So, we assume that $G$ is defined over 
$\BF_p$, of split type, with Frobenius map $F$. We assume that 
$T \subset B$ are both $F$-stable, and that
$L \subset P$ are $F$-stable.  
Thus $F$ acts trivially on $W$ and on $W_L$. 
We first note that
\begin{lem} 
Let $a \in N_W(W_L)$ and choose $\da \in N_G(T) \cap N_G(L)$.
Assume that $\da \in Z_G(u)$.  Then $\ad \da$ stabilizes $\CB_u^L$, and
acts on $H^*(\CB_u^L)$ in such a way that 
$\ad \da (w)  = awa\iv$ for $w \in W_L$.
\end{lem}
\begin{proof}
Since $\da \in N_G(L)$, $\da$ acts on $\CB^L$ by the adjoint action 
$\ad \da$, which
stabilizes $\CB_u^L$ since $\da \in Z_G(u)$.  Hence $\da$ acts 
naturally on $H^*(\CB_u^L)$.  In order to compare this action with
the action of $W_L$, we shall recall the construction of Springer 
representations of $W_L$.  Let 
\begin{equation*}
\wt L = \{ (x, gB) \in L \times \CB^L \mid g\iv xg \in B\},
\end{equation*}
and $\pi: \wt L \to L$ be the first projection.  Let $L_r$ be 
the set of regular semisimple elements in $L$.  Then $\pi\iv(L_r)$
is isomorphic to 
\begin{equation*}
 \wt L_r = T_r \times L/T, 
\end{equation*}
where $T_r = T \cap L_r$. 
Let $\pi_0: \wt L_r \to L_r$ be the map defined by 
$\pi_0 : (t, gT) \mapsto g\iv tg$, which coincides with the 
restriction of $\pi$ on $\wt L_r$ under the identification 
$\pi\iv(L_r) \simeq \wt L_r$. 
Then $\pi_0$ is an unramified Galois covering with 
group $W_L$, and for a constant sheaf $\Ql$ on $\wt L_r$, 
$\CL = \pi_*\Ql$ is a $W_L$-equivariant local system 
on $L_r$.  Thus $K = \IC(L, \CL)$ is a $W_L$-equivariant complex 
on $L$, and it is known by Lusztig that $K \simeq \pi_*\Ql$.  Thus 
 for each $u \in L$, the stalk $\CH^i_u(K)$ at $u$ of the
$i$-th cohomology sheaf of $K$ gives rise to a 
$W_L$-module $H^i(\CB_u^L)$.
\par
Now $\da$ acts on $\wt L_r$ (resp. on $L_r$) by 
$\ad \da: (t, gT) \mapsto (\da t\da\iv, \da g\da\iv T)$ 
(resp. $\ad \da: x \mapsto \da x\da\iv$), and $\pi_0$ commutes with 
$\ad \da$.  Hence $\CL$ becomes 
an $\da$-equivariant local system.  Since 
$\pi_0\iv(t) = \{ (wtw\iv, wT) \mid w \in W_L\}$ for $t \in T_r$,
the stalk $\CL_t$ has a natural structure of the regular $W_L$-module.  
Then the isomorphism $\CL_{\da t\da\iv} \to \CL_t$ is given by 
$\ad \da\iv$ under the identification $\CL_{x} \simeq \Ql[W_L]$ for 
$x \in L_r$.  It follows that $\CL$ is 
$\lp \da\rp\ltimes W_L$-equivariant, where $\lp\da\rp$ is 
a cyclic group generated by $\da$, and $\da$ acts on $W_L$ by 
$\ad \da(w) = awa\iv$. 
By the functoriality of $\IC$ functor, $K$ turns out 
to be a $\da$-equivariant complex on $L$ under the adjoint action of 
$\da$, which is regarded as a 
$\lp \da\rp \ltimes W_L$-equivariant complex on $L$.
Hence for $u \in L$ such that $\da u\da\iv =u$, $\CH^i_u(K)$ has a 
structure of $\lp \da\rp \ltimes W_L$-module. 
\par
On the other hand, $\da$ acts naturally on $\wt L$ and on $L$ 
by the adjoint action, which commute with $\pi$. 
Thus $\pi_*\Ql$ is $\da$-equivariant, which is isomorphic to $K$ 
as the complex with $\da$-action. Hence the action of $\da$
on $\CH^i_u(K)$ coincides with the action on $H^i(\CB^L_u)$ 
induced from the adjoint action of $\da$ on $\CB_u^L$.
The lemma follows from this.
\end{proof}
Next we show the following lemma.
\begin{lem} 
There exists a representative $\da \in N_G(T)\cap N_G(L) \cap Z_G(u)$ 
such that 
$\da$ acts trivially on $H^*(\CB_u)$ 
and that $\da^e$ acts trivially on $H^*(\CB_u^L)$.
In particular, $H^*(\CB_u^L)$ has a structure of $\wt W_L$-module. 
\end{lem}
\begin{proof}
First consider the case (a) in 1.6. 
 Let $H$ be the subgroup of $G$ generated 
by $U_{\a}$ with $\a \in \F_{L'}$, where $U_{\a}$ is the root subgroup
corresponding to $\a$.  Then $H$ is a connected reductive subgroup of
$L'$ whose Weyl group coincides with $W_{L'}$.  Since 
$H \subset Z_G(u)$, we have $H \subset Z_G^0(u)$.  
One can choose a representative $\da \in N_H(T_1)$ of $a \in W_{L'}$, 
where $T_1$ is a maximal torus of $H$ contained in $T$. 
Then $\da \in Z_G^0(u) \cap N_G(L)$ and $\da^e \in  T_1$. 
 Since 
$T_1 \subset Z_G(u)$, we see that $T_1 \subset Z_L^0(u)$.  Thus, 
$\da^e \in Z_L^0(u)$.  Hence $\da$ satisfies the condition.
\par
Next consider the case (b) in 1.6. 
Let $L_1$ be the Levi subgroup containing $L$ of type 
$X_{n_0} + A_{en_1-1} + \cdots + A_{en_r-1}$.
We have a natural projection $\pi: L_1 \to \bar L_1 = L_1/Z^0(L_1)$, and 
an isogeny map 
$\th : \wt L_1 = G_0 \times SL_{en_1}\times\cdots\times SL_{en_r} 
           \to \bar L_1$, 
where $G_0$ is the simply connected semisimple group of type 
$X_0$. 
Put $\bar u = \pi(u) \in \bar L_1$.   
Now $Z_{L_1}(u)$ acts on $H^*(\CB_u)$.  Since $Z^0(L_1)$ acts 
trivially on $H^*(\CB_u)$, we have an action of 
$Z_{L_1}(u)/Z^0(L_1) = Z_{\bar L_1}(\bar u)$ on $H^*(\CB_u)$.
Let $\wt u$ be an element in $\wt L_1$ such that $\th(\wt u) = \bar u$.
$\wt u = (u_0, u_1, \dots, u_r)$ can be chosen as given in 1.4.
We choose $\dda \in \wt L_1$ as follows;  
put $\dda = (a_0, a_1, \dots, a_r)$ with $a_0 \in G_0$, and 
$a_i \in SL_{en_i}$ for $1 \le i \le r$.  We put $a_0 = 1$ and
choose $a_1, \dots, a_r$ so that $a_i \in Z^0_{SL_{en_i}}(u_i)$ 
and that $a_i^e \in Z(SL_{en_i})$.  Such a choice is always possible
for $u_i$ of type $(n_i, \dots, n_i)$.  Thus 
$\dda \in Z^0_{\wt L_1}(\wt u)$.  It follows that $\th(\dda)$ is contained 
in a connected subgroup of $Z_{\bar L_1}(\bar u_1)$, and by the previous
remark, $\th(\dda)$ acts trivially on $H^*(\CB_u)$. 
Now take $\da \in Z_{L_1}(u)$ such that $\pi(\da) = \th(\dda)$.  Then 
$\da \in N_G(T) \cap N_G(L)$, and acts trivially on $H^*(\CB_u)$.
On the other hand, similar to $\pi,\th$, we have a map 
$\pi': L \to \bar L = L/Z^0(L)$ and 
$\th': \wt L = G_0 \times (SL_{n_1})^e 
        \times\cdots \times (SL_{n_r})^e \to \bar L$.
Let $\bar u = \pi'(u) \in \bar L$, and $\wt u \in \wt L$ such that 
$\bar u = \th'(\wt u)$.  Then we have an isomorphism 
$H^*(\CB^L_u) \simeq H^*(\CB_{\bar u}^{\bar L}) 
        \simeq H^*(\CB_{\wt u}^{\wt L})$ 
compatible with the actions of $Z_L(u), Z_{\bar L}(\bar u)$ and 
$Z_{\wt L}(\wt u)$ with respect to $\pi',\th'$.
We have $\dda^e \in Z(SL_{n_1})^e\times Z(SL_{n_2})^e \times \cdots$. 
Since  the action of 
$Z(SL_{n_1})^e \times Z(SL_{n_2})^e \times\cdots$ can be extended to 
an action of $Z(GL_{n_1})^e \times Z(GL_{n_2})^e\times\cdots$ on 
$H^*(\CB_{\wt u}^{\wt L})$, $\dda^e$ acts trivially on 
$H^*(\CB_{\wt u}^{\wt L})$, and so $\da^e$ acts trivially on 
$H^*(\CB_u^L)$.
\end{proof}
\para{2.4.}
Let $\CZ = Z_L^0$ be the identity component of the center of $L$.
Put $\CB_{\CZ} = \{ B' \in \CB \mid \CZ \subset B' \}$.  Then 
$\CB_{\CZ}$ is decomposed into connected components  
\begin{equation*}
\CB_{\CZ} = \coprod_{d \in W_L\backslash W}\CB_{\CZ, d}, 
\end{equation*}
where $\CB_{\CZ,d} = \{ {}^{xd}B \mid x \in L\}$, which is 
isomorphic to $\CB^L$ under the map $B' \mapsto B' \cap L$.
\par
Put 
\begin{equation*}
\CZ\reg = \{ z \in \CZ \mid Z^0_G(z) = L \}.
\end{equation*}
Then for any $t \in \CZ\reg$, we have $\CB_t = \CB_{\CZ}$ by 
Lemma 2.2 (c) in [L3], and so 
$\CB_{tu} = \CB_u \cap \CB_t = \CB_u \cap \CB_{\CZ}$.
It follows that 
\begin{equation*}
\CB_{tu} = \coprod_{d \in W_L\backslash W}(\CB_{\CZ,d}\cap \CB_u),
\end{equation*} 
where $\CB_{\CZ,d} \cap \CB_u$ is isomorphic to $\CB_u^{L}$ under the 
map $B' \mapsto B' \cap L$.  This implies that 
\begin{equation*}
\tag{2.4.1}
H^{2n}(\CB_{tu}) \simeq 
       \bigoplus_{d\iv \in W/W_L} H^{2n}(\CB_{\CZ,d} \cap \CB_u).
\end{equation*}
The right hand side of (2.4.1) has a natural structure 
of the induced $W$-module $\Ind_{W_L}^WH^{2n}(\CB_u^L)$.  
It is proved in [L3, Proposition 1.4] that (2.4.1) is actually 
an isomorphism of $W$-modules. 
Let $a \in W$ be as in the theorem.  
Since $\da \in N_G(L)$, it 
stabilizes $\CZ$, and so $\da$ acts on $\CB_{\CZ}$ via $\ad\da$.  
It is easy to see that $\da$ induces a permutation action 
on the components of 
$\CB_{\CZ}$; $\da : \CB_{\CZ, d} \mapsto \CB_{\CZ, ad}$.  
It follows that $\da$ induces an automorphism on  
$H^{2n}(\CB_{tu})$, which maps the factor 
corresponding to $d\iv \in W/W_L$ to $d\iv a\iv \in W/W_L$.
Under the isomorphism 
$H^{2n}(\CB_{\CZ,d} \cap \CB_u) \simeq H^{2n}(\CB_u^L)$, 
the factor corresponding to $d\iv \in W/W_L$ is written as
$d\iv\otimes H^{2n}(\CB_u^L)$, and $\da$ maps 
$d\iv\otimes H^{2n}(\CB_u^L) \to 
            d\iv a\iv \otimes H^{2n}(\CB_u^L)$.
On the other hand, by Lemma 2.3,  $\da^e$ acts trivially on 
$H^{2n}(\CB_u^L)$, and induces an action of $\wt W_L$ on it. 
Hence $\da$ induces an action of $\vG_0$ on 
$H^{2n}(\CB_{tu}) \simeq \Ind_{W_L}^W H^{2n}(\CB_u^L)$, which is 
given by 
$\da: d\iv \otimes x \mapsto d\iv a\iv \otimes \da x$ for each factor
$d\iv\otimes H^{2n}(\CB_u^L)$. 
\par
Now we define an action of $\vG$ on $H^*(\CB_{tu})$ by 
$a: x \mapsto \z^n\da x$ for $x \in H^{2n}(\CB_{tu})$, where
$\da x$ is the action of $\vG_0$ on $H^{2n}(\CB_{tu})$ 
given as above.  Since the action
of $\da \in G$ commutes with that of $W$, $H^*(\CB_{tu})$ turns out 
to be a $\vG \times W$-module, which we denote by 
$H^*(\CB_{tu})^{[\z]}$.
The following lemma is immediate from the above discussion.
\begin{lem} 
There exists an isomorphism of $\vG \times W$-modules
\begin{equation*}
H^*(\CB_{tu})^{[\z]} \simeq \gind_{W_L}^WH^*(\CB_u^L)^{(\z)}.
\end{equation*}
\end{lem}
In view of Lemma 2.5, in order to prove the theorem it is enough 
to show the following proposition.
\begin{prop} 
Under an appropriate choice of (a good prime) $p$,  
there exists an isomorphism 
of $\vG \times W$-modules for any $t \in \CZ_r$,
\begin{equation*}
H^*(\CB_u)^{(\z)} \simeq H^*(\CB_{tu})^{[\z]}.
\end{equation*}
\end{prop}
\para{2.7.}
The remainder of this section is devoted to the proof of the proposition.
We shall prove it by modifying the arguments in [L3].
By [Sh1], [Sh2], [BS], the following 
fact is known; assume that $G$ is simple modulo center.  
Then for each unipotent class $C$ of $G$, there 
exists $u_1 \in C^F$, called a split unipotent element,
 such that $F$ acts on $H^{2n}(\CB_{u_1})$ as a
scalar multiplication by $p^n$. (In the case where $G$ is of type
$E_8$, we assume that $p \equiv 1 \pmod 4$). 
 Since the component group
$A_G(u_1) = Z_G(u_1)/Z_G^0(u_1)$ is isomorphic to $S_3, S_4, S_5$ or 
$(\BZ/2\BZ)^k$ for some $k$, there exists a positive integer 
$s_0$ (independent of $p$) such that $F^{s_0}$ acts on 
$H^{2n}(\CB_u)$ by a scalar multiplication by $p^{s_0n}$ for 
any unipotent element $u$ of $G^F$ (e.g., one can take $s_0 = |S_5|$.)  
Similarly, $F^{s_0}$ acts on $H^{2n}(\CB^L_u)$ by a scalar 
multiplication by $p^{s_0n}$ for any unipotent element $u \in L^F$.  
Note that the isomorphism in (2.4.1) is $F$-equivariant. 
Hence $F^{s_0}$ acts also as a scalar
multiplication by $p^{s_0n}$ for $H^{2n}(\CB_{tu})$. 
\par
Note that $\da$ acts trivially on $H^{2n}(\CB_u)$ by Lemma 2.3.  
It follows that one can write 
\begin{align*}
\tag{2.7.1}
\Tr((F^s\da)^i w, H^*(\CB_u)) &= \sum_{n \ge 0}a_n(w)p^{isn},  \\
\tag{2.7.2}
\Tr((w,a^i), H^*(\CB_u)^{(\z)}) &= \sum_{n \ge 0}a_n(w)\z^{in},  
\end{align*}
for any $w \in W, 0 \le i \le e-1$ and for any positive integer $s$ divisible
by $s_0$, where $a_n(w) = \Tr(w, H^n(\CB_u))$ 
are integers for each $n \ge 0$. 
\par
On the other hand, by the description of the action of $F$ and of
$\da$ on $H^n(\CB_{tu})$ in 2.4, together with Lemma 2.5, one can
write
\begin{align*}
\tag{2.7.3}
\Tr((F^s\da)^i w, H^*(\CB_{tu})) &= \sum_{n \ge 0}b_{n,i}(w)p^{isn}, \\
\tag{2.7.4}
\Tr((w,a^i), H^*(\CB_{tu})^{[\z]}) &= \sum_{n \ge 0}b_{n,i}(w)\z^{in},
\end{align*}
for $w, i, s$ as above, where $b_{n,i}(w)$ are 
certain integers.
\par
For an integer $x$ and a prime number $l$, 
we denote by $m_{l}(x)$ the multiplicative
order of $x$ in $\BZ/l\BZ$, i.e., the smallest positive integer
$m$ such that $x^m \equiv 1 \pmod l$. 
The following is a key for the proof of Proposition 2.6.
\begin{lem} 
Assume that $p \equiv 1 \pmod 4$. Let $s_0, e$ be fixed positive 
integers coprime to $p$. 
Then there exist infinitely many prime numbers $l$ satisfying the
following properties.
\begin{enumerate}
\item
$m_{l}(p^s) = e$ for a certain integer $s$ divisible by $s_0$.
\item
$l-1$ is divisible by $e$.
\end{enumerate}
\end{lem}
\begin{proof}
By our assumption, the image of $s_0e$ on $\BF_p = \BZ/p\BZ$ is
non-zero.  Hence the map $x \mapsto s_0ex + 1$ induces a bijective map 
on $\BF_p$.  Thus there exists $c \in \BZ$ such that
the image of $s_0ec + 1$ in $\BF_p$ is  contained 
in $\BF_p^* - (\BF_p^*)^2$.  Put $\a = s_0ec +1$.  
Then $\a$ is prime to $p$, and so 
$(\a-1)p$ and $\a$ are coprime each other.  Then by Dirichlet's
theorem on arithmetic progression, there exist infinitely many 
prime numbers $l$ of the form $l = n(\a-1)p +\a$ for some positive 
integer $n$.  
It is enough to show that these $l \ge 3$ satisfy the assertion 
of the lemma.
For an integer $a$ and a prime number $p$,
let $\displaystyle\bigl(\frac{a}{p}\bigr)$ be the 
Legendre symbol, i.e.,
\begin{equation*}
\bigl(\frac{a}{p}\bigr) = \begin{cases}
             1 &\quad\text{ if $x^2 \equiv a\pmod p$ for some $x \in \BZ$}, \\ 
             -1 & \quad\text{ otherwise.}
                          \end{cases}
\end{equation*}
We show that 
\begin{equation*}
\tag{2.8.1}
\bigl(\frac{p}{l}\bigr) = -1.
\end{equation*}
In fact, by the quadratic reciprocity law (e.g., [Se]), 
we have
\begin{equation*}
\bigl(\frac{p}{l}\bigr)\bigl(\frac{l}{p}\bigr)
   = (-1)^{\frac{p-1}{2}\cdot\frac{l-1}{2}} = 1.
\end{equation*} 
The second equality follows from the assumption that 
$p \equiv 1 \pmod 4$. 
Hence we have 
$\displaystyle\bigl(\frac{p}{l}\bigr) 
        = \bigl(\frac{l}{p}\bigr)$.
But $l \equiv \a \pmod p$, and so 
$\displaystyle\bigl(\frac{l}{p}\bigr) 
   = \bigl(\frac{\a}{p}\bigr) = -1$
since the image of $\a$ is not contained in $\BF_p^2$ by 
our choice of $\a$.  Hence (2.8.1) holds. 
\par
Now (2.8.1) is equivalent to $p^{(l-1)/2} \equiv -1 \pmod l$.
It follows that $m_{l}(p) = l-1$. 
Since $l-1 = s_0ec(np+1)$, we see that $m_{l}(p^s) = e$ for
$s = s_0c(np+1)$ and that $l-1$ is divisible by $e$.  
Thus this $l$ satisfies the assertion of the 
lemma.  The lemma is proved. 
\end{proof}
\para{2.9}
For  given integers $s_0 \ge 1, e \ge 2$, we choose 
a prime number $p$ such that $p$ is not a factor of $e, s_0$ and that
$p \equiv 1 \pmod 4$, and fix it once and for all.  
For a multiple $s$ of $s_0$, put $F' = F^s\da$ and $q = p^s$. 
Under the setting in 1.6, we
shall describe the set $\CZ\reg$ more precisely. 
As in [L3, Lemma 2.2], $\CZ\reg$ can be written as 
$\CZ\reg = \CZ - \bigcup_{\b}\ker(\b|_{\CZ})$, where 
$\b$ runs over all the roots in $\F - \F_L$. 
($\b|_{\CZ}$ gives a non-trivial character of $\CZ$
for $\b \in \F - \F_L$).
\par
First consider the case (a).
Let $L'\der$ be the derived subgroup of $L'$, and $S'$ be the
split maximal torus of $L'\der$ contained in $T$.  Then 
$S' \subset \CZ$.  Put $S'\reg = S' \cap \CZ\reg$.
Now $W_{L'}$ leaves the set $\F - \F_L$ invariant. 
For each $\b \in \F - \F_L$, put 
$H_{\b} = \bigcap_{x \in \vG}\ker(x(\b)|_{S'})$.  Then 
$H_{\b}$ is an $F'$-stable subgroup of $S'$, and we see that 
\begin{equation*}
\tag{2.9.1}
{S'}\reg^{F'} = {S'}^{F'} - \bigcup_{\b \in \F - \F_L}H_{\b}^{F'}.
\end{equation*}
$H_{\b}$ is a closed subgroup of $S'$, and we put 
$e_{\b} = |H_{\b}/H_{\b}^0|$ for each $\b \in \F - \F_L$.
\par
Let $\CP'$ be the set of all prime numbers $l$ satisfying the 
condition in Lemma 2.8.  Thus $\CP'$ is an infinite set. 
We denote by $\CP$ the subset of $\CP'$ consisting of $l$ such that
$l > |\F - \F_L|$ and that $l$ does not divide $e_{\b}$ 
($\b \in \F - \F_L) $.
Thus $\CP$ is an infinite set also. 
\par
Next we consider the case (b).  
We may assume that $G$ has a connected center of dimension 1, and 
that the derived subgroup of $G$ is simply connected, almost simple. 
Let $k$ be an algebraic closure of $\BF_q$.
We see that there exists a subtorus $S$ of $\CZ$ such that 
$S \simeq (k^*)^c$, where
$c$ is the number of irreducible components of $\F_L$.
Since $a$ permutes the factors $k^*$ in $S$, we see
that ${S}^{F'} \simeq (\BF^*_{q^e})^r\times (\BF_q^*)^{r'}$, 
where $r'$ is equal to 1 or 0 according to the cases where 
$X_0$ is non-empty or empty.  
Since $\vG \subset N_W(W_L)$, $\vG$
preserves the set $\F - \F_L$.  
For each $\b \in \F - \F_L$, put 
$K_{\b} = \bigcap_{x \in \vG}\ker (x(\b)|_{S})$.
Then $K_{\b}$ is an $F'$-stable subgroup of $S$, and we have
\begin{equation*}
\tag{2.9.2}
{S}^{F'}\reg = {S}^{F'} - \bigcup_{\b \in \F - \F_L}K_{\b}^{F'},
\end{equation*}
where ${S}\reg = S \cap \CZ\reg$.
$K_{\b}$ is a closed subgroup of $S$, and put 
$e_{\b} = |K_{\b}/K^0_{\b}|$ for each $\b \in \F - \F_L$.
Under the identification 
$S^{F'} \simeq (\BF_{q^e}^*)^r \times (\BF_q^*)^{r'}$, 
we see that 
$K_{\b}^{0F'} \simeq (\BF^*_{q^e})^{r-1} \times (\BF_q^*)^{r'}$ or 
$K_{\b}^{0F'} \simeq \BF^*_{q^{e'}}\times 
  (\BF^*_{q^e})^{r-1} \times (\BF_q^*)^{r'}$, 
where $e'$ is a proper divisor of $e$.
 (Let $S_i$ be the subtorus of $S$ corresponding to the factor
$eA_{n_i-1}$ for $i = 1, \dots, r$.  Then the former case occurs 
if $\b|_{S_i}, \b|_{S_j}$ are non-trivial for some $i \ne j$, 
and the latter case occurs if $\b|_{S_i}$ is non-trivial 
for only one $i$.  Note that by our assumption in 1.4, $\b$
is non-trivial on $S_1\times\cdots\times S_r$.)
\par
Let $\CP'$ be as in the case (a).  We define a subset $\CP$ of $\CP'$ as 
the set of prime numbers $l \in \CP'$ such that $l > |\F - \F_L|$
and that $l$ does not divide $e_{\b}$.
\par 
The next lemma is a variant of Lemma 3.4 in [L3].
\begin{lem} 
Assume that $l \in \CP$, and let $s$ be a multiple of $s_0$ such that 
$m_l(p^s) = e$ (see Lemma 2.8).  Put $F' = F^s\da$.  Then 
there exists $t \in \CZ\reg$ such that $F'(t) = t$ and that 
$t^l = 1$.
\end{lem}
\begin{proof}
First consider the case (a) in 1.6.
It is enough to show, for each $l \in \CP$, 
that there exists $t \in {S'}^{F'}\reg$ such that $t^l = 1$.  
Note that $a$ is a regular element of order $e$ in $W_{L'}$.  Put 
$V = \BR\otimes_{\BZ}X(S')$.  Thus $W_{L'}$ acts on $V$ as a 
reflection group.  Let $\z$ be a primitive $e$-th root of unity, 
and let $a(e)$ be the dimension of the 
eigenspace $V(a,\z) \subset V$ of $a$ with eigenvalue $\z$.
We show that 
\begin{equation*}
\tag{2.10.1}
\sharp\{ t \in {S'}^{F'} \mid t^l = 1\} = l^{a(e)}.
\end{equation*}
\par
By a general formula, we have 
$|{S'}^{F'}| = |\det_V(qI - a)| = P_a(q)$, where $P_a(x)$ is the
characteristic polynomial of $a \in W_{L'}$. 
Since $a$ is regular $P_a(x)$ can be written, by [Sp1, 4.2], as 
\begin{equation*}
P_a(x) = \F_e(x)^{a(e)}\F'(x),
\end{equation*}
where $\F_e(x)$ is the cyclotomic polynomial of degree $e$, and 
$\F'(x)$ is a product of cyclotomic polynomials $\F_{e'}(x)$ with 
$e' < e$.
By our assumption $m_l(q) = e$, $\F_e(q)$ is divisible by $l$, 
and $\F'(q)$ is not divisible by $l$. 
This means that each minimal $F'$-stable torus $M$ of $S'$ corresponding 
to the factor $\F_{e}(x)$ contains an element of order $l$.
Since $\{ t \in M^{F'} \mid t^l = 1\} \subset \BF^*_{q^e}$, 
$M^{F'}$ contains exactly $l$
elements $t$ such that $t^l = 1$.
Thus (2.10.1) is proved. 
\par
For $\b \in \F - \F_L$, let $V_{\b}$ be the subspace of $V$ 
which is orthogonal to $x(\b)$ for all $x \in \vG$.  Then $V_{\b}$ 
can be identified with $\BR\otimes_{\BZ} X(H_{\b}^0)$.   
$\vG$ stabilizes $V_{\b}$, and let $V_{\b}(a, \z)$ be the eigenspace
of $a$ on $V_{\b}$ with eigenvalue $\z$.  Since $a$ is $L$-regular, 
we have 
$\dim V_{\b}(a, \z) < \dim V(a, \z) = a(e)$.
It follows that the characteristic polynomial $P'_a(x)$ of $a$ on 
$V_{\b}$ contains the factor $\F_e(x)$ with multiplicity less than
$a(e)$.
By a similar argument as above,  minimal $F'$-stable subtori of
$H^0_{\b}$ corresponding to $\F_e(x)$ only contain elements of order $l$. 
This implies that 
\begin{equation*}
\sharp\{ t \in H_{\a}^{F'} \mid t^l = 1\} = 
    \sharp\{ t \in H_{\a}^{0F'} \mid t^l = 1\} \le l^{a(e) -1}.
\end{equation*}
It follows, by (2.9.1), that 
\begin{align*} 
\sharp\{ t \in {S'}\reg^{F'} \mid t^l = 1\} &= 
       \sharp\{ t \in {S'}^{F'} \mid t^l = 1, 
      t \notin \bigcup_{\b \in \F - \F_L} H_{\b}^{F'} \}  \\
  & \ge l^{a(e)} - Nl^{a(e) - 1} = l^{a(e) - 1}(l - N), 
\end{align*}
where $N = |\F - \F_L|$.
Since $l > N$ by our assumption, there exists $t \in {S'}\reg^{F'}$ such 
that $t^l = 1$.  This proves the lemma in the case (a).
\par
Next consider the case (b) in 1.6.  
It is enough to show, for each $l \in \CP$,  that there exists 
$t \in S\reg^{F'}$ such that $t^l = 1$. 
We note that 
$q^{e'}-1$ is not divisible by $l$ for any divisor $e' < e$ of $e$
by the assumption $m_l(q) = e$.
Since $S^{F'} \simeq (\BF^*_{q^e})^r \times (\BF_q^*)^{r'}$
(cf. 2.9), we have 
\begin{equation*}
\sharp\{ t \in S^{F'} \mid t^l = 1\} = l^r.
\end{equation*}
We consider $K_{\b}$ given in 2.9.  By   
the discussion in 2.9, we have
\begin{equation*}
\sharp\{ t \in K^{F'}_{\b} \mid t^l = 1\} 
= \sharp\{ t \in K^{0F'}_{\b} \mid t^l = 1\} = l^{r-1}.
\end{equation*}
It follows, by (2.9.2), that 
\begin{align*}
\sharp\{ t \in S^{F'}\reg \mid t^l = 1\}
    &= \sharp\{ t \in S^{F'} \mid t^l = 1, 
          t \notin \bigcup_{\b \in \F - \F_L}K_{\b}^{F'}\} \\
    &\ge l^r - Nl^{r-1} = l^{r-1}(l - N),
\end{align*}  
where $N$ is as before. Since $l > N$ by our assumption,  
the lemma holds also for the case (b).
\end{proof}
We need the following lemma due to Lusztig.
\begin{lem}[{[L3, Lemma 3.2]}]
Let $H$ be a finite group, and $\f$ a virtual character of 
$H$ (over a field of characteristic 0).  
Assume that $\f$ is integral valued.  Let $x, y \in H$
be such that $xy = yx$ and $y^l = 1$ for a prime number $l$.
Then $\f(xy) - \f(x) \in l\BZ$. 
\end{lem}
\para{2.12.}
Let $s_0$ be as in 2.7, and $\CP$ be as in 2.9.
Let $F' = F^s\da$ be as in Lemma 2.10 for a fixed 
$l \in \CP$.  Let $R_{w, i} = R_{T_w}(1)$ be the 
Deligne-Lusztig's virtual character of $G^{F'^i}$ for 
$i = 1, \dots, e$, where $T_w$ 
is an $F'^i$-stable maximal torus of $G$ corresponding to 
$w \in W \simeq W(T_1)$ (here $W(T_1) = N_G(T_1)/T_1$ for an 
$F'$-stable pair $T_1 \subset B_1$). 
Let us choose $t \in \CZ\reg$ as in Lemma 2.10.
Then we have 
\begin{align*}
\tag{2.12.1}
\Tr(F'^iw, H^*(\CB_u)) &= \Tr(u, R_{w,i}), \\
\Tr(F'^iw, H^*(\CB_{tu})) &= \Tr(tu, R_{w,i}).
\end{align*}
We remark that (2.12.1) was proved in [L2] under the assumption that 
$p^s$ is large enough (which is determined only by the data of the Dynkin 
diagram of $G$).  Thus if we replace $s_0$ in 2.7 by a suitable large
number, the result in [L2] is applicable.  One can also apply 
[Sh3, Theorem 2.2] instead of [L2], where the restriction on $p^s$ is removed.
\par
Since $R_{w,i}$ are integral valued, one can apply Lemma 2.11 for 
$H = G^{F'^i}$ and $x = u, y = t$. Hence we have
\begin{equation*}
\Tr(u, R_{w,i}) = \Tr(tu, R_{w,i}) \mod l\BZ.
\end{equation*}
It follows from (2.12.1) that
\begin{equation*}
\tag{2.12.2}
\Tr(F'^iw, H^*(\CB_u)) = \Tr(F'^iw, H^*(\CB_{tu}))
                           \mod l\BZ.
\end{equation*}
\par
Let $\z_0$ be a fixed primitive $e$-th root of unity in $\BC$, and
$R$ the ring of integers of the cyclotomic field $\BQ(\z_0)$.
Let $\CI$ be the set of non-zero prime ideals $\Fp$ in $R$ such that 
$\Fp$ contains one of the numbers $1 - \z_0^i$ for $i = 1, \dots, e-1$
and $\z_0$.
Let $\bar\CI$ be the set of prime numbers $l$ such that 
$\Fp \cap \BZ = l\BZ$ for $\Fp \in \CI$.  
Since $\CI$ is a finite set, $\bar\CI$ is a finite set.
So, $\CP - \bar\CI$ is an infinite set.  
Let $\CJ$ be the set of prime ideals $\Fp$ of $R$ such that 
$\Fp \cap \BZ = l\BZ$ with $l \in \CP - \bar\CI$.
Then $\CJ$ is an infinite set. 
Now $R/\Fp$ is a finite extension of $\BF_l$.
Let $\bar\z_0$ be the image of $\z_0$ in $R/\Fp$.
Since $l \in \CP$, $l-1$ is divisible by $e$.  Hence 
$\bar\z_0 \in \BF_l^*$, which has order $e$ by our choice of $\Fp$.
Since $m_l(p^s) = e$, the image of $p^s$ in $\BZ/l\BZ$ has order $e$. 
Hence there exists $j$ such that 
\begin{equation*}
\tag{2.12.3}
p^s - \z_0^j \in \Fp.
\end{equation*}
Note that the number $j$ is determined by the choice of $\Fp$, which 
we denote by $j(\Fp)$.
For $j = 1, \dots, e-1$, let $\CJ_j$ be the set of prime ideals 
$\Fp$ in $\CJ$ such that $j(\Fp) = j$.  Thus $\CJ = \bigcup_j \CJ_j$,
and so there exists $j_0$ such that $\CJ_0 = \CJ_{j_0}$ is an infinite set.
We put $\z = \z_0^{j_0}$.  By (2.12.3), $\z$ is a primitive $e$-th 
root of unity.  
\par
We remark that $H^*(\CB_{tu}) = H^*(\CB_u \cap \CB_{\CZ})$ is 
independent of the choice of $t \in \CZ\reg$.
Then in view of (2.7.1) $\sim$ (2.7.4), together with (2.12.3), we see that 
\begin{align*}
\Tr((F^s\da)^iw, H^*(\CB_u)) &= 
        \Tr((w,a^i), H^*(\CB_u)^{(\z)}) \mod \Fp,  \\ 
\Tr((F^s\da)^iw, H^*(\CB_u \cap \CB_{\CZ})) &= 
                     \Tr((w,a^i), H^*(\CB_u\cap \CB_{\CZ})^{[\z]}) \mod \Fp
\end{align*}
for any $\Fp \in \CJ_0$.
Combined with (2.12.2), we have
\begin{equation*}
\Tr((w,a^i), H^*(\CB_u)^{(\z)}) 
     = \Tr((w,a^i), H^*(\CB_u\cap \CB_{\CZ})^{[\z]}) \mod \Fp
\end{equation*}
for $\Fp \in \CJ_0$.
Since $\CJ_0$ is an infinite set, we conclude that 
\begin{equation*}
\Tr((w,a^i), H^*(\CB_u)^{(\z)}) 
     = \Tr((w,a^i), H^*(\CB_u\cap \CB_{\CZ})^{[\z]}).
\end{equation*}
Hence Proposition 2.6 is proved, and the theorem follows. 
\par\bigskip  
\section{Applications}
\para{3.1.}
Let $W_L$ be the subgroup of $W$, and $\vG$ the 
subgroup of $W$ generated by $a \in N_W(W_L)$ such that 
$\vG$ and $W_L$ generate the semidirect product group 
$\wt W_L = \vG \ltimes W_L$. 
Let $V = V^{(\z)}$ be the $\vG \times \wt W_L$-module as in 1.3. 
(We write $\vG$ as $\vG_0$ if it is regarded 
as a subgroup of $\wt W_L$,  cf. 1.3.)
Then 
$V$ can be decomposed as 
$V = \bigoplus_{i \in \BZ/e\BZ}V^{(i)}$, where $V^{(i)}$ is the eigenspace of 
$a \in \vG$ with eigenvalue $\z^i$, which is a $\wt W_L$-submodule of $V$.
Then we have
\begin{align*}
\Ind_{W_L}^WV &= \bigoplus_i
      \bigoplus_{w \in W/\wt W_L}\bigoplus_j
                   wa^j\otimes V^{(i)}  \\
              &= \bigoplus_i\bigoplus_{w \in W/\wt W_L}
                    \bigoplus_k wb_k\otimes V^{(i)},
\end{align*}
where $b_k = \sum_j\z^{jk}a^j \in \BC[\vG]$ (the group ring of $\vG$).
For each $i \in \BZ$, let $\p^{(i)}$ the linear character of $\vG$ defined by 
$\p^{(i)}(a) = \z^i$.  Then $\vG$-module $\BC b_k$ is afforded by 
$\p^{(-k)}$.  Let $V^{(i)}_n$ be the eigenspace of $a \in \vG_0$ on 
the $\wt W_L$-module $V^{(i)}$  
with eienvalue $\z^n$.    
Let $(\gind_{W_L}^W V)^{(k)}$ be the eigensapce of $a \in \vG$ with eigenvalue
$\z^k$.  Then we have the following lemma. 
\begin{lem}  
\begin{enumerate}
\item  Let the notations be as above.  We have
\begin{equation*}
\tag{3.2.1}
(\gind_{W_L}^WV)^{(k)} \simeq \bigoplus_{w \in W/\wt W_L}
              \bigoplus_{j \in \BZ/e\BZ}\bigoplus_{0 \le n < e}
                   wb_{k-n-j}\otimes V_n^{(j)}
\end{equation*}
as vector spaces.  In particular, 
$\dim (\gind_{W_L}^WV^{(\z)})^{(k)}$ is independent of the 
choice of $k \in \BZ/e\BZ$, which is given by  
\begin{equation*}
\tag{3.2.2}
\dim (\gind_{W_L}^WV)^{(k)} = [W: \wt W_L]\dim V.
\end{equation*}
\item
Assume that $\vG$ commutes with $W_L$.  Then we have 
\begin{equation*}
\tag{3.2.3}
(\gind_{W_L}^WV)^{(k)} \simeq 
           \bigoplus_{ j \in \BZ/e\BZ}
                 \Ind_{\vG\times W_L}^W(\p^{(-k+j)}\otimes V^{(j)}) 
\end{equation*}
as $W$-modules.  
\end{enumerate}
\end{lem}
\begin{proof}
Under the action of $\vG$ on 
$\gind_{W_L}^W V$, $wb_k\otimes V^{(i)}_n$ is contained in an 
eigenspace of $a$ with eigenvalue $\z^{k+n+j}$.
Then (i) follows easily from the discussion in 3.1.
Now assume that $\vG$ commutes with $W_L$.  Then 
$b_k\otimes V^{(i)}$ has a structure of $\vG \times W_L$-module
given by $\p^{(-k)}\otimes V^{(i)}$.  (ii) follows from the formula
(3.2.1) by noticing that $V^{(i)} = V_0^{(i)}$.  The lemma is proved.
\end{proof} 
We consider a Levi subgroup $L \subset G$ and a unipotent element
$u \in L$, and take $\vG = \lp a\rp \subset N_W(W_L)$ satisfying the 
condition in 1.6.
We apply the preceding argument to the situation 
$V^{(\z)} = H^*(\CB_u^L)^{(\z)}$.    
Then as a corollary to Theorem 1.7, we have
\begin{prop} 
Under the setting in Theorem 1.7, we have, for $0 \le k \le e-1$,  
\begin{equation*}
\tag{3.3.1}
\bigoplus_{n \equiv k \mod e}H^{2n}(\CB_u) \simeq 
                    \bigoplus_{w \in W/\wt W_L}
              \bigoplus_{j \in \BZ/e\BZ}\bigoplus_{0 \le n < e}
                   wb_{k-n-j}\otimes 
                      H^{2j}(\CB_u^L)_n
\end{equation*} 
as vector spaces, where $b_i \in \BC[\vG]$ and $H^{2n}(\CB_u^L)_n$
is the eigenspace of $a \in \vG_0$ with eigenvalue $\z^n$. 
In particular, 
$\dim \bigl(
      \bigoplus_{n \equiv k \mod e}H^{2n}(\CB_u)\bigr)$ is 
independent of the choice of $k$.
In the case \rm{(a)} in 1.6, (3.3.1) can be made more precise as follows;
\begin{equation*}
\tag{3.3.2}
\bigoplus_{n \equiv k \mod e}H^{2n}(\CB_u) \simeq 
                    \Ind_{\vG \times W_L}^W \biggl(
              \bigoplus_{j \in \BZ/e\BZ}
                   \p^{(-k+j)}\otimes 
                      H^{2j}(\CB_u^L)\biggr)
\end{equation*} 
as $W$-modules.  
\end{prop}
\begin{proof}
By Theorem 1.7  $\gind_{W_L}^WH^*(\CB_u^L)^{(\z)}$ is 
isomorphic to $H^*(\CB_u)^{(\z)}$ as $\vG \times W$-modules.
Since $(H^*(\CB_u)^{(\z)})^{(k)} 
        = \bigoplus_{n \equiv k \mod e}H^{2n}(\CB_u)$,
the corollary follows from Lemma 3.2.
\end{proof}
\remarks{3.4.}
(i) \ 
In the case where $u = 1$, the cohomology ring 
$H^*(\CB_u) = H^*(\CB)$ coincides 
with the coinvariant algebra $R$ of $W$.
In the special case where $G$ is of type $A_{n-1}$, i.e.,
$W \simeq \FS_n$,
we consider $W_L \simeq \FS_{n-re}$ for $1 \le e \le n$. Then
$W_{L'}\simeq \FS_{re}$, and if we choose a regular element
$a \in W_{L'}$ as a product 
of disjoint cycles of length $e$,  Proposition 3.3 can be applied. 
This recovers the formula obtained by Morita and Nakajima [MN1]. 
\par
More generally, consider the Weyl group $W$ acting on the real vector space 
$V$ as the reflection module.  For $v \in V$, let $W_v$ be the 
stabilizer of $v$ in $W$, and $N_v$ the stabilizer of the line 
$\BR v$ in $W$.  Note that $W_v$ is normal in $N_v$, and 
$W_v$ coincides with $W_L$ for a certain Levi subgroup in $G$.  
Then for any $\vG = \lp a\rp$ such that 
$\vG \subset N_v$, Bonnaf\'e, Lehrer and Michel [BLM] have proved
a similar formula as in Proposition 3.3.  So our formula (3.3.2) 
can be regarded as a special case of theirs. (Note that they treat 
a more general case, where $W$ is a complex reflection group and $\vG$ 
is not necessarily cyclic, in a framework of coinvariant algebras.)
\par
(ii) \ 
We consider a unipotent element $u \in L$ in the case where
$G = GL_n$. $u \in G$ can  be written as $u = u_{\mu}$ 
by a partition $\mu = (1^{m_1}, 2^{m_2}, \dots, n^{m_n})$ of $n$.
Take a positive integer $e \ge 2$, and let 
$I$ be a subset of $\{ 1, \dots, n\}$ 
such that $e \le m_i$ for $i \in I$.
We consider a Levi subgroup $L$ of 
type $X_0 + e\sum_{i \in I}A_{i-1}$, where 
$X_0 = A_k$ with $k = \sum_{i \notin I} im_i + \sum_{i \in I}i(m_i-e) -1$. 
Then we have $W_{L'} = \{1\}$, and one can choose $u \in L$ so that 
it satisfies the assumption of the case (b) in 1.6. 
Thus Proposition 3.3 can be applied.  This covers the results on 
the stability of dimensions obtained in [MN2], [MN3], where they considered 
the case $|I| = 1$ or
the case all the $m_i$ are divisible by $e$. 
\par\medskip
Returning to the general setup, we consider the case 
where $u$ is a regular unipotent element in 
$L$.  Then $H^*(\CB_u^L) = H^0(\CB_u^L) \simeq \BC$ is a trivial
$W_L$-module.  Thus Proposition 3.3 implies the following.
\begin{cor}  
Let $G$ be a simple algebraic group modulo center, and $L$ a Levi
subgroup in $G$.  Let $u$ be a regular unipotent element in $L$.
Let $\vG = \lp a\rp$ be a subgroup of $N_W(W_L)$ of order $e$ 
satisfying the conditions in 1.6.
Then for $k = 0, \dots, e-1$, we have
\begin{equation*}
\bigoplus_{n \equiv k \mod e}H^{2n}(\CB_u) \simeq
          \Ind_{\vG\ltimes W_L}^W \wt\p^{(-k)}
\end{equation*}
as $W$-modules, where $\wt\p^{(-k)}$ is the character of 
$\vG \ltimes W_L$ obtained as the pull back of $\p^{(-k)}$
under the projection $\vG \ltimes W_L \to \vG$.
\end{cor}
\begin{proof}
In the case (a), the assertion follows from (3.3.2).  So we consider 
the case (b).  In the setup of 3.1, $V^{(i)}$ is a trivial $W_L$-module
$\BC$ for $i = 0$ and zero otherwise. Then we see that 
$V^{(0)} = V^{(0)}_0$, and $b_k \otimes V^{(0)}$ has 
a structure of $\wt W_L$-module $\wt\p^{(-k)}$. The assertion 
follows from the formula in 3.1. 
\end{proof}
 \para{3.6.}
Let $G$ be a simple algebraic group defined over $\BF_q$ 
with Frobenius map $F$.  We assume that $G^F$ 
is of split type.  
The Green function $Q_{T_w}$ is defined as the restriction of 
the Deligne-Lusztig's virtual  character $R_{T_w}(1)$ to 
the set of unipotent elements in $G^F$.  
We assume that $p = \ch \BF_q$ is good, and in the case where $G$ is of type
$E_8$, we further assume that $q \equiv 1 \pmod 4$.
Then as explained in 2.7, 
for each unipotent class $C$ of $G$, there exists a split element 
$u \in C^F$.  As in 2.12, we have 
\begin{equation*}
\tag{3.6.1}
Q_{T_w}(u) = \sum_{n \ge 0}\Tr(w, H^{2n}(\CB_u))q^n.
\end{equation*} 
Hence there exists a polynomial $\BQ_{w,C}(x) \in \BZ[x]$ such that 
$Q_{T_w}(u) = \BQ_{w,C}(q)$. 
Concerning the values of Green functions at root of unity, we have
the following.
\begin{prop}  
Suppose that $G, L$ and $u \in L$ are as in Corollary 3.5.  Then we have
\begin{equation*}
\tag{3.7.1}
\BQ_{w,C}(\z^j) = |W_L|\iv
        \sharp\{ x \in W \mid x\iv wx \in a^jW_L\}
\end{equation*}
for $j = 0, \dots, e-1$.  
In particular, the value $\BQ_{w,C}(\z')$ is independent of the 
choice of a primitive $e$-th root of unity $\z'$.
\end{prop} 
\begin{proof}
Put $c_i(w) = \sharp\{ x \in W \mid x\iv wx \in a^iW_L\}$ for 
$i = 0, \dots, e-1$.
Then 
\begin{align*}
\bigl(\Ind_{\vG\ltimes W_L}^W\wt\p^{(-k)}\bigr)(w)
    &= |\vG\ltimes W_L|\iv\sum_{i=0}^{e-1}
                 \sum_{\substack{x \in W \\ x\iv wx \in a^iW_L}}
                        \wt\p^{(-k)}(x\iv wx)  \\
    &= |\vG\ltimes W_L|\iv \sum_{i=0}^{e-1}c_i(w)\z^{-ki}.
\end{align*}
It follows, by (3.6.1) together with Corollary 3.5, that 
\begin{align*}
\BQ_{w,C}(\z^j) &= 
    \sum_{k=0}^{e-1}\z^{kj}\sum_{n \equiv k \mod e}\Tr(w, H^{2n}(\CB_u))  \\ 
     &= |\vG\ltimes W_L|\iv \sum_{i=0}^{e-1}c_i(w)  
                         \sum_{k=0}^{e-1}\z^{(j-i)k}  \\  
                   &= |W_L|\iv c_j(w).
\end{align*}
Hence we obtain the formula (3.7.1). 
Let $\z^j$ be a primitive $e$-th root of unity.  There exists
an element $\t \in \Gal(\BQ(\z)/\BQ)$ such that $\t(\z) = \z^j$.
By (3.6.1), we see that $\BQ_{w,C}(\z) \in \BQ(\z)$ and that 
$\t(\BQ_{w,C}(\z)) = \BQ_{w,C}(\z^j)$.  But since 
$\BQ_{w,C}(\z) \in \BZ$ by (3.7.1), we conclude that 
$\BQ_{w,C}(\z) = \BQ_{w,C}(\z^j)$.  
This proves the proposition.
\end{proof}
\remark{3.8.}
In the case where $G = GL_n$ and $L$ is of type 
$A_{m-1} + \cdots + A_{m-1}$ ($e$-times) with $n = em$, take a regular
unipotent element $u$ in $L$.  Then $u = u_{\mu} \in G$ with $\mu = (m^e)$. 
For $w \in W = \FS_n$, let $\la(w) = (1^{l_1}, 2^{l_2}, \dots)$ be 
the partition of $n$ corresponding to the cycle decomposition of
$w$. 
Then one can show by a direct computation (cf. [M, (6.2)])
that 
\begin{equation*}
|W_L|\iv\sharp\{ x \in W \mid x\iv wx \in aW_L\}
   = \begin{cases}
           e^{l(\la(w))}  &\quad
        \text{ if $e \mid l_i$ for all $i$}, \\
           0            &\quad\text{ otherwise,}
     \end{cases} 
\end{equation*} 
where $l(\la)$ is the number of parts for a partition $\la$.
Thus we recover the formula in [LLT, Theorem 3.2, Theorem 3.4] 
concerning the values of Green polynomials of $GL_n$ at roots of unity.
\para{3.9.}
We give some more examples where Proposition 3.3 can be applied.
\par
(i) \ Assume that $G$ is of type $B_n$ and $L$ is a Levi subgroup of
type $B_m$ with $m < n$.  Then $L'$ is of type $A_{n-m-1}$.  
For any $u \in L$ and a divisor $e$ of $n-m$, the 
proposition can be applied. Similar results hold
also for $C_n$ or $D_n$.
\par
(ii) \ Assume that $G = Sp_{2n}$.  Then a unipotent element $u \in G$
can be written as $u = u_{\mu}$ as an element of $GL_{2n}$, where
$\mu = (1^{m_1}, 2^{m_2}, \dots)$ is a partition of $2n$ such that 
$m_i$ is even for odd $i$. 
Take an even integer $e \ge 2$, and let $I$ be a subset of 
odd integers $\{ 1, 3, \dots, 2n-1\}$ such that $e \le m_i$ for $i \in I$.
We consider a Levi subgroup $L$ of type 
$X_0 + e\sum_{i \in I}A_{i-1}$, where $X_0$ is of type $C_k$ 
with $2k = \sum_{i\notin I}im_i + \sum_{i \in I}i(m_i-e)$.   
Then $W_{L'} = \{ 1\}$, and as in Remarks 3.4 (ii), one can find
$u \in L$ so that the case (b) in 1.6 can be applied. 
Similar results hold also for type $B_n$ and $D_n$.
\par
(iii) \ Assume that $G$ is of type $E_7$, and choose $L$ of type 
$A_2$ so that $L'$ is of type $A_4$.  Take any unipotent element 
$u \in L$.  Then the proposition can be applied with $e = 5$. 
\par
\end{document}